# Discrete Darboux Transformation for Discrete polynomials of hypergeometric Type

Gaspard Bangerezako

September 26, 1997

*Abstract.* Darboux Transformation, well known in second order differential operator theory, is applied here to the difference equations satisfied by the discrete hypergeometric polynomials(Charlier, Meixner-Kravchuk, Hahn).

*Keywords* : Discrete Darboux Transformation, Orthogonal Polynomials.

## 1. Introduction

Since Darboux, who showed how $z = Ay + By'$ solves $z'' = (\psi + h)z$ when $y$ satisfies $y'' = (\varphi + h)y$ [1], numerous generalizations have been investigated.
We consider here a second order difference equation

$$H(x;j)\Phi(x;j) = 0 \qquad (1)$$

where
$$H(x;j) = E^2 + v(x;j)E + u(x;j) \qquad (2)$$

with
$$E\Phi(x;j) = \Phi(x+1;j) , \qquad (3)$$

$x \in R,\ j \in Z.$

Suppose that one can form the products

$$\begin{aligned} H(x;j) - \mu(j) &= (E + g(x;j))(E + f(x;j)) \\ H(x;j+1) - \mu(j) &= (E + f(x;j))(E + g(x;j)) + \alpha(j) \end{aligned} \qquad (4)$$

then, the operator $H(x;j+1)$ is called *Discrete Darboux Transformation* of $H(x;j)$. From (4), we have the following commutation relation

$$H(x;j+1)(E+f(x;j)) = (E+f(x;j))(H(x;j)+\alpha(j)) \qquad (5)$$

which is a discrete analog of the so called *dressing chain* [5],[7].
The dressing chain (5) is equivalent to the system

$$\begin{aligned} f(x;j)+g(x+1;j) &= f(x+1;j+1)+g(x;j+1) \\ f(x;j)\,g(x;j)+\alpha(j) &= f(x;j+1)\,g(x;j+1)+\mu(j+1)-\mu(j) \end{aligned} \qquad (6)$$

In the continuous case, many of the questions concerning intrinsic structure (Hamiltonian, integrability,..) of such chains were elucided in [5],[4]. In the discrete case, similar structures remain generally obscure. Even some particular considerations of such systems appearing today are mainly directed to the cases when the shift operator in (3) acts on $j$ but not on $x$. This is typically the case when one is treating the discrete Schrödinger problem or particularly the polynomial recurrence relations. Note that except some particular cases of selfadjointness of $x$ and $j$, this nuance is very significant. The polynomial recurrence relations case for example, is characterized by the linearity of the eigenvalue, which is then $x$, which a priori facilitates the application of the Darboux transformation techniques.

In the next section, we shall prove that some finite difference hypergeometric type operators are particular solutions of the chain(5) . This leads to difference and recurrence relations for the corresponding eigenfunctions.[1]

In the last section, we shall apply our result to the classical orthogonal polynomials of a discrete variable on a linear lattice (Charlier, Meixner-Kravchuk, Hahn). The Charlier and Meixner-Kravchuk cases having been treated similarly in [2], we firstly succeed to handle the Hahn case, specialized by the nonlinearity of the eigenvalue, as a function of $j$ .

## 2. Finite Difference Analogs of $\lambda(n)$-Eigenfunctions of Hypergeometric Type

Let $\Phi(x;n)$ , $x \in R$ , $n \in Z$ , be a given system of hypergeometric

---
[1] For another type of factorization construction, see [6].

functions such that

$$(\sigma(x)\Delta\nabla + \tau(x)\Delta)\,\Phi(x;n) = \lambda(n)\,\Phi(x;n) \qquad (7)$$

where $\sigma(x)$ and $\tau(x)$ are polynomials of degree less or equal than 2 and 1 respectively, $\lambda(n) = n\tau' + \frac{1}{2}n(n-1)\sigma''$ [3].
For convenience, we shall call them $\lambda(n)$-*eigenfunctions*. It is clear that this set includes discrete polynomials of hypergeometric type. Let's note that in this case, the number $n$ is the degree of the corresponding polynomial.

One can easily find some function $\rho(x)$ such that (7) is equivalent to

$$(E^2 - [2\sigma(x+1) + \tau(x+1) + \lambda(n)]E + (\sigma(x) + \tau(x))\sigma(x+1))(\rho(x)\Phi(x;n)) = 0 \qquad (8)$$

with

$$\frac{\rho(x+1)}{\rho(x)} = \sigma(x) + \tau(x).$$

Let $L = E^2 - [2\sigma(x+1) + \tau(x+1)]E + (\sigma(x) + \tau(x))\sigma(x+1)$ and

$$H(x;n) = L - \lambda(n)E = E^2 - [2\sigma(x+1) + \tau(x+1) + \lambda(n)]E + (\sigma(x) + \tau(x))\sigma(x+1)$$

Supposing the existence of two polynomials $f(x;n)$ and $g(x;n)$ of second degree with identical leading coefficients such that

$$H(x;n) - \mu(n) = (E + g(x;n))(E + f(x;n)) \qquad (9)$$

for some constant $\mu(n)$, one can verify that

$$(E + f(x;n))(E + g(x;n)) = H(x;n') - \mu(n) \qquad (10)$$

where

$$\lambda(n') = \lambda(n) + \Delta(f(x;n) - g(x;n)), \qquad (11)$$

$n'$ being some function of $n$, which will be determined latter.
(9) and (10) give

$$H(x;n')(E + f(x;n)) = (E + f(x;n))H(x;n) \qquad (12)$$

In order to determine $f(x;n)$ and $g(x,n)$, one needs to note that equation (9) leads to the system

$$\begin{aligned}
f(x+1;n) + g(x;n) &= -2\sigma(x+1) - \tau(x+1) - \lambda(n) \\
f(x;n)g(x;n) &= (\sigma(x) + \tau(x))\sigma(x+1) - \mu(n)
\end{aligned} \qquad (13)$$

which is in fact a discrete Riccati equation.

Setting

$$f(x;n) = -\sigma(x) - \tau(x) - \frac{1}{2}\lambda(n) + \varphi(x;n)$$
$$g(x;n) = -\sigma(x+1) - \frac{1}{2}\lambda(n) - \varphi(x+1;n) \qquad (14)$$

the first equation in (13) will be automatically verified. The second reads

$$\frac{1}{2}\lambda(n)(\sigma(x+1) + \sigma(x) + \tau(x)) + \frac{1}{4}\lambda^2(n) + \mu(n) + (\sigma(x) + \tau(x))\varphi(x+1;n)$$
$$-\sigma(x+1)\varphi(x;n) + \frac{1}{2}\lambda(n)\Delta\varphi(x;n) - \varphi(x;n)\varphi(x+1;n) = 0 \qquad (15)$$

a discrete Riccati equation in relation to $\varphi(x;n)$. Looking for polynomial solutions of degree $\leq 1$, $\varphi(x;n) = \phi(n)x + \psi(n)$ ; knowing that $\sigma(x) = \sigma_0 x^2 + \sigma_1 x + \sigma_2$ , $\tau(x) = \tau_0 x + \tau_1$ and equating coefficients in left side of (15) to zero, one finds two possible sets of solutions

$$\phi_1(n) = \tau_0 + (n-1)\sigma_0 \;;\; \phi_2(n) = -n\sigma_0 \qquad (16)$$

$$\psi_{1,2}(n) = \frac{\phi_{1,2}(n)(\tau_1 + \tau_0 - \sigma_0 - \phi_{1,2}(n)) + \lambda(n)\sigma_0 + \lambda(n)\sigma_1 + \frac{1}{2}\lambda(n)\tau_0}{2\phi_{1,2}(n) + 2\sigma_0 - \tau_0} \qquad (17)$$

$$\mu_{1,2}(n) = \psi_{1,2}(n)(\psi_{1,2}(n) + \phi_{1,2}(n) + \sigma_1 + \sigma_0 - \tau_1) - \frac{1}{2}\lambda(n)(\sigma_0 + \sigma_1 + 2\sigma_2 + \tau_1)$$
$$-\phi_{1,2}(n)(\sigma_2 + \tau_1 + \frac{1}{2}\lambda(n)) - \frac{1}{4}\lambda^2(n) \qquad (18)$$

On the other side, (14) reads

$$f(x;n) = -\sigma_0 x^2 + (\phi(n) - \sigma_1 - \tau_0)x + \psi(n) - \sigma_2 - \tau_1 - \frac{1}{2}\lambda(n) \qquad (19)$$

$$g(x;n) = -\sigma_0 x^2 - (\phi(n) + 2\sigma_0 + \sigma_1)x - \sigma_0 - \sigma_1 - \sigma_2 - \frac{1}{2}\lambda(n) - \phi(n) - \psi(n) \qquad (20)$$

Thus, the conditions advanced in (9) are all satisfied .
Next, from (19) , (20), we get $\Delta(f-g) = 2\phi(n) + 2\sigma_0 - \tau_0$, and using (16), it follows ($\lambda(n) = n\tau_0 + n(n-1)\sigma_0$),

$$(\Delta(f-g))_1 = 2\phi_1(n) + 2\sigma_0 - \tau_0 = \tau_0 + 2n\sigma_0 = \lambda(n+1) - \lambda(n) \qquad (21)$$

$$(\Delta(f-g))_2 = 2\phi_2(n) + 2\sigma_0 - \tau_0 = -(\tau_0 + 2(n-1)\sigma_0) = \lambda(n-1) - \lambda(n) \qquad (22)$$

Refering to (11), this means that we have proved that $n'_{1,2} = n \pm 1$ and (12) reads
$$H(x; n \pm 1)(E + f_{1,2}(x; n)) = (E + f_{1,2}(x; n))H(x; n) \qquad (23)$$
which is the searched commutation relation(5)$(j := n)$.
From (23) and (8), it obviously follows that for any $\lambda(n)$-eigenfunction $\Phi(x; n)$ of hypergeometric type, the following difference relations are valid
$$\begin{aligned} c_1(n)\,\tilde\Phi(x; n+1) &= (E + f_1(x; n))\tilde\Phi(x; n) \\ c_2(n)\,\tilde\Phi(x; n-1) &= (E + f_2(x; n))\tilde\Phi(x; n) \end{aligned} \qquad (24)$$
where $\tilde\Phi(x; n) = \rho(x)\Phi(x; n)$,
$$f_{1,2}(x; n) = -\sigma_0 x^2 + (\phi_{1,2}(n) - \sigma_1 - \tau_0)x + \psi_{1,2}(n) - \sigma_2 - \tau_1 - \tfrac{1}{2}\lambda(n).$$

From this, of course, recurrence relations, for $\Phi(x; n)$, can be deduced.

3. **Examples**

$$Charlier$$

$$H(x; n) = E^2 - (x + \mu + \lambda(n) + 1)E + \mu(x+1) \;;\; \rho(x) = \mu^x \;;\; \lambda(n) = -n$$

$$\begin{aligned} f_1(x; n) &= -x + n & ; & f_2(x; n) &= -\mu \\ g_1(x; n) &= -\mu & ; & g_2(x; n) &= -x + n - 1 \\ \mu_1(n) &= \mu n + \mu & ; & \mu_2(n) &= \mu n \end{aligned}$$

$$Meixner$$

$$\begin{aligned} H(x; n) &= E^2 - [(\mu+1)x + \mu(\gamma+1) + 1 + \lambda(n)]E + \mu x^2 + \mu(\gamma+1)x + \gamma\mu \;; \\ & \rho(x) = \mu^x \Gamma(x+\gamma) \;;\; \lambda(n) = -n(1-\mu) \end{aligned}$$

$$\begin{aligned} f_1(x; n) &= -x + n & ; & f_2(x; n) &= -\mu(x + \gamma + n) \\ g_1(x;) &= -\mu(x + \gamma + n - 1) & ; & g_2(x; n) &= -x + n - 1 \\ \mu_1(n) &= \mu(n\gamma + n^2 + n + \gamma) & ; & \mu_2(n) &= \mu n(\gamma + n - 1) \end{aligned}$$

$$Hahn$$

$$\begin{aligned}
H(x;n) &= E^2 + [2x^2 + (6 + \beta - \alpha - 2N)x + (5 + 2\beta - \alpha - 3N - \beta N - \lambda(n))]E + \\
&\quad [x^4 + (4 + \beta - \alpha - 2N)x^3 + (6 + 3\beta - 3\alpha - 6N + N^2 - 2N\beta + \alpha N - \alpha\beta)x^2 + \\
&\quad (4 + 3\beta - 3\alpha - 6N + 2N^2 - 4N\beta + 2N\alpha - 2\alpha\beta + N^2\beta + N\alpha\beta)x + \\
&\quad 1 + \beta - \alpha - 2N + N^2 - 2N\beta + \alpha N - \alpha\beta + N^2\beta + N\alpha\beta] ; \\
&\rho(x) = \tfrac{\Gamma(x+\beta+1)}{\Gamma(-x+N)} ; \quad \lambda(n) = -n(n + \alpha + \beta + 1)
\end{aligned}$$

$$\begin{aligned}
f_1(x;n) &= x^2 - (N + \alpha + n - 1)x - (\beta+1)(N-1) - \tfrac{1}{2}\lambda(n) + \psi_1(n) \\
g_1(x;n) &= x^2 + (3 + n + \beta - N)x + 2 + \beta + n - N - \tfrac{1}{2}\lambda(n) - \psi_1(n) \\
\mu_1(n) &= \psi_1(n)(\psi_1(n) - 1 - \beta N - n) - \tfrac{1}{2}\lambda(n)(\beta+1)(N-1) - \tfrac{1}{2}\lambda(n)(N + \alpha - 1) \\
&\quad + (n + \alpha + \beta + 1)(\beta + 1)(N - 1) + \tfrac{1}{4}\lambda(n)(n+2)(n + \alpha + \beta + 1)
\end{aligned}$$

$$\begin{aligned}
f_2(x;n) &= x^2 + (2 + \beta - N + n)x - (\beta+1)(N-1) - \tfrac{1}{2}\lambda(n) + \psi_2(n) \\
g_2(x;n) &= x^2 + (2 - n - N - \alpha)x - N - \alpha - n + 1 - \tfrac{1}{2}\lambda(n) - \psi_2(n) \\
\mu_2(n) &= \psi_2(n)(\psi_2(n) + n + \alpha - \beta N + \beta) - \tfrac{1}{2}\lambda(n)(N + \alpha - 1) - n(\beta+1)(N-1) \\
&\quad - \tfrac{1}{2}\lambda(n)(\beta+1)(N-1) - \tfrac{1}{2}\lambda(n)n - \tfrac{1}{4}\lambda^2(n)
\end{aligned}$$

$$\begin{aligned}
\psi_1(n) &= \tfrac{(n+\alpha+\beta+1)(\beta+1)(N-1) - \lambda(n)(N+\alpha) + \tfrac{1}{2}\lambda(n)(\alpha+\beta+2)}{2 + 2n + \alpha + \beta} \\
\psi_2(n) &= \tfrac{n(\beta+1)(N-1) + \lambda(n)(N+\alpha) - \tfrac{1}{2}\lambda(n)(\alpha+\beta+2)}{2n + \alpha + \beta}
\end{aligned}$$

It is clear that the same technique can be applied as well to the q-versions of the preceding polynomials. Extension to other difference operators is in progress.

Aknowledgements: We thank Prof. A. P. Magnus for having suggested the present explorations and for fruitful discussions. Thanks are addressed as well to the Belgium General Agency for Cooperation with Developing Countries(AGCD) for financial support.

*Permanent address* : Université du Burundi, Faculté des sciences, Département de Mathématique, B. P. 2700 Bujumbura, Burundi.

*Current address* : Institut de Mathématique, Université Catholique de Louvain, Chemin du Cyclotron, 2 , 1348 Louvain-La-Neuve, Belgium.

*E-mail address* : bangerezako@agel.ucl.ac.be